\documentclass{siamart251216}

\usepackage{lipsum}
\usepackage{amsfonts}
\usepackage{graphicx}
\usepackage{epstopdf}
\usepackage{algorithmic}
\usepackage{booktabs}
\usepackage{amsopn}

\DeclareMathOperator*{\argmin}{arg\,min}
\newcommand{\vect}[1]{\mathbf{#1}} 
\newcommand{\mat}[1]{\mathbf{#1}} 
\newcommand{\norm}[1]{\|#1\|}
\usepackage[most]{tcolorbox}
\usepackage{subcaption}
\usepackage{caption}
\pdfoutput=1

\usepackage[colorinlistoftodos,prependcaption,textsize=tiny]{todonotes}
\ifpdf
  \DeclareGraphicsExtensions{.eps,.pdf,.png,.jpg}
\else
  \DeclareGraphicsExtensions{.eps}
\fi

\usepackage{datetime}
\newdateformat{monthyeardate}{%
  \monthname[\THEMONTH] \THEDAY, \THEYEAR}

\makeatletter
\newcommand*{\addFileDependency}[1]{
  \typeout{(#1)}
  \@addtofilelist{#1}
  \IfFileExists{#1}{}{\typeout{No file #1.}}
}
\makeatother

\usepackage{academicons}
\usepackage{xcolor}

\usepackage{enumitem}
\setlist[enumerate]{leftmargin=.5in}
\setlist[itemize]{leftmargin=.5in}

\newsiamthm{problem}{Problem}
\newsiamthm{assumption}{Assumption}
\newsiamremark{remark}{Remark}
\newsiamremark{hypothesis}{Hypothesis} 
\crefname{hypothesis}{Hypothesis}{Hypotheses}
\newsiamremark{example}{Example}
\newsiamthm{claim}{Claim}
\newsiamthm{conjecture}{Conjecture}

\headers{Cumulative Edge-Detection for Image Reconstruction}{
Okunola T., Pasha M., Kilmer M.
}

\title{An Efficient Cumulative Edge-Detection Method for Image Reconstruction
}
\author{
Toluwani Okunola\thanks{Department of Mathematics, Tufts University} (\email{toluwani.okunola@tufts.edu})
\and 
Mirjeta Pasha\thanks{Department of Mathematics, Virginia Tech} (\email{mpasha@vt.edu})
\and
Misha E. Kilmer\thanks{Department of Mathematics, Tufts University}    (\email{misha.kilmer@tufts.edu})
}
\usepackage{amsopn}

\newtcolorbox{keybox}[1]{
  colback=blue!5!white,
  colframe=black,
  title=#1,
  fonttitle=\bfseries,
  boxrule=0.5pt,
  arc=2mm
}

\begin{document}

\maketitle

\begin{abstract}
When reconstructing images from noisy measurements, such as in medical scans or scientific imaging, we face an inverse problem: recovering an unknown image from indirect, corrupted observations. These problems are typically ill-posed, meaning small amounts of noise can lead to inaccurate reconstructions. Regularization techniques address this by incorporating prior assumptions about the solution, such as smoothness or sparsity. However, standard methods often blur sharp edges—the boundaries between tissues or structures—losing critical detail.

A powerful strategy for edge preservation is iterative reweighting, which solves a sequence of weighted subproblems with adaptively updated weights. Non-cumulative schemes derive weights from the current iterate alone and can be solved efficiently using the Recycled Majorization-Minimization Generalized Krylov Subspace method (RMM-GKS). The cumulative approach of Gazzola et al. progressively accumulates edge information across iterations, achieving superior edge preservation but at high computational cost.

This work introduces CR-$\ell_q$-RMM-GKS, which combines cumulative edge detection with computational efficiency. We integrate Gazzola's cumulative weighting with RMM-GKS, which handles general $\ell_q$ penalties ($0 < q \le 2$), automatically selects regularization parameters, and recycles Krylov subspaces between iterations, reducing the nested structure to two levels.

Numerical experiments in signal deblurring and tomography demonstrate that CR-$\ell_q$-RMM-GKS produces significantly sharper edge reconstructions than standard non-cumulative methods. In particular, CR-$\ell_1$-RMM-GKS outperforms both standard $\ell_1$ methods and CR-$\ell_2$-RMM-GKS, indicating that cumulative weighting and $\ell_1$ penalties are highly complementary.
\end{abstract}
\begin{keywords}
edge preservation, Krylov subspace methods, iterative regularization, inverse problems.
\end{keywords}

\begin{AMS}
65F22, 65F10, 65K10, 68U10.
\end{AMS}

\section{Introduction} \label{sec:intro}
In many scientific and medical imaging applications—from computed tomography (CT) to astronomical imaging—we aim to reconstruct an underlying image $\vect{x}_{\text{true}} \in \mathbb{R}^{n}$ from indirect, noisy measurements $\vect{b} \in \mathbb{R}^{m}$ \cite{hansen2010discrete}. The linear forward problem is
\begin{equation}
\vect{b} = \mat{A}\vect{x}_{\text{true}} + \boldsymbol{\eta},
\label{eq:forward}
\end{equation}
where $\mathbf{A} \in \mathbb{R}^{m \times n}$ (with $m \geq n$ in the overdetermined setting; 
underdetermined cases such as sparse-angle tomography are also 
accommodated), represents the forward imaging model and $\boldsymbol{\eta} \in \mathbb{R}^{m}$ is additive noise. A key challenge is that $\mat{A}$ is typically severely ill-conditioned, with singular values that decay rapidly toward zero. This means that small perturbations in the data $\vect{b}$—even tiny amounts of noise—can cause massive changes in the solution. Consequently, the naive least-squares solution, $\argmin_{\vect{x}} \|\mat{A}\vect{x} - \vect{b}\|_2^2$, amplifies noise dramatically and produces unusable reconstructions.

To obtain a stable and meaningful solution, we must employ \emph{regularization}, which incorporates prior knowledge or assumptions about the desired solution $\vect{x}_{\text{true}}$. One of the most well-known regularization methods is Tikhonov regularization, which solves
\begin{equation}
\vect{x}_{\lambda} = \argmin_{\vect{x}} \left\{ \|\mat{A}\vect{x} - \vect{b}\|_2^2 + \lambda \|\mat{L}\vect{x}\|_2^2 \right\}.
\label{eq:tikhonov}
\end{equation}
Here, $\lambda > 0$ is the regularization parameter balancing data fidelity and regularization, and $\mat{L}$ is typically a discrete gradient or derivative operator. The regularization term $\|\mat{L}\vect{x}\|_2^2$ penalizes large values in $\mat{L}\vect{x}$, which correspond to sharp changes in pixel intensities. This promotes smoothness throughout the image, effectively suppressing noise in regions that should be uniform. 

\begin{figure}[htbp]
    \centering
    \includegraphics[width=0.9\linewidth]{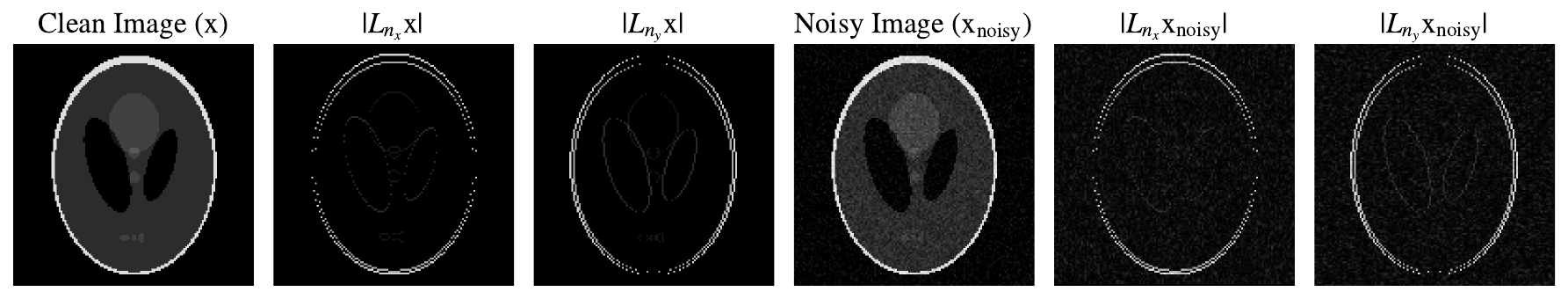}
    \caption{Visualizing the effect of the discrete gradient operator $\mat{L}$. From left to right: a clean image, its horizontal and vertical gradient magnitudes, a noisy version, and its gradient magnitudes. Both actual edges and noise contribute to large gradient values, complicating edge preservation in standard $\ell_2$ regularization.}
    \label{fig:description_of_Lx}
\end{figure}

However, sharp intensity changes are also exactly what occur at \emph{edges}—the boundaries between different structures, tissues, or materials that carry crucial diagnostic or scientific information. As illustrated in Figure~\ref{fig:description_of_Lx}, both true edges and noise produce large gradient values. Standard $\ell_2$ regularization cannot distinguish between them and therefore applies uniform smoothing indiscriminately, blurring critical structural details. {Preserving edges while suppressing noise is therefore a central challenge in image reconstruction.}

\subsection{Edge-Preserving Regularization via Iterative Reweighting}

To overcome the edge-blurring limitation of standard Tikhonov regularization, a powerful general strategy has emerged: \emph{iterative reweighting}. The core idea is to solve a sequence of weighted $\ell_2$ subproblems, where the weights are updated to adapt the penalty to the structure of the evolving solution. By assigning smaller weights at detected edge locations and larger weights elsewhere, these methods can preserve edges while still regularizing smooth regions.

Iterative reweighting schemes differ primarily in \emph{how the weights are computed}, leading to two major classes:

\paragraph{Non-cumulative weighting} 
The most widely used approach derives weights from the \emph{current iterate alone} \cite{bardsley2012, candes2008, daubechies2010, wohlberg2007}. A particularly successful example uses $\ell_q$-norm penalties ($0 < q < 2$), most commonly $q=1$ as in Total Variation (TV) regularization~\cite{rudin1992}. Since the $\ell_1$-norm is non-smooth, Majorization-Minimization (MM) methods~\cite{hunter2004, lange2004, rodriguez2009} convert it into a sequence of weighted $\ell_2$ subproblems. Modern solvers such as the Recycled Majorization-Minimization Generalized Krylov Subspace (RMM-GKS) method~\cite{pasha2023} achieve impressive computational efficiency by recycling Krylov subspaces across iterations and automatically selecting regularization parameters. Operator-splitting approaches such as Split Bregman / 
ADMM~\cite{goldstein2009, boyd2011} also solve TV-type problems 
efficiently but require a different algorithmic infrastructure 
(dual variable updates, parameter tuning for the augmented 
Lagrangian) and do not naturally accommodate the Krylov 
recycling strategy central to our framework.

\paragraph{Cumulative weighting}
A fundamentally different philosophy, introduced by Gazzola et al.~\cite{gazzola2020}, updates weights by \emph{accumulating} edge information across iterations. Rather than discarding the edge detection history, their method builds weights that progressively ``remember'' previously identified edges through a multiplicative update mechanism. This cumulative strategy excels at preserving fine structures and produces superior edge quality compared to non-cumulative methods. However, the method requires solving multiple nested optimization problems at each step, making it computationally expensive for large-scale problems.

\subsection{Contribution}

This work bridges the gap between computational efficiency and edge preservation quality by introducing \emph{Cumulative Reweighted $\ell_q$ Regularization with Recycled Majorization-Minimization Generalized Krylov Subspace} (CR-$\ell_q$-RMM-GKS). Our framework combines the cumulative edge-detection strategy of Gazzola et al.~\cite{gazzola2020} with the computational efficiency of the RMM-GKS solver~\cite{pasha2023}, which handles general $\ell_q$-norm penalties ($0 < q \le 2$), automatically selects regularization parameters, and recycles Krylov subspaces between outer iterations.  In this respect, our approach is related to hybrid Krylov 
methods~\cite{chung2017} that combine iterative projection with automatic 
parameter selection, but extends this paradigm to cumulative $\ell_q$ 
reweighting with subspace recycling across outer iterations.

At each outer iteration, we update cumulative weights and pass the modified regularization operator to RMM-GKS, which warm-starts from the Krylov subspace constructed in the previous iteration. This reduces the deeply nested structure of prior cumulative methods to two levels while extending cumulative weighting to general $\ell_q$ penalties. The resulting method solves
\[
\vect{x}^{(\ell)}
=
\argmin_{\vect{x}}
\left\{
\|\mat{A}\vect{x}-\vect{b}\|_2^2
+
\lambda_\ell
\|\mat{D}^{(\ell)}\mat{L}\vect{x}\|_q^q
\right\},
\qquad 0<q\le2,
\]
where each subproblem benefits from both automatic parameter selection and subspace recycling.

Numerical experiments on 1D deblurring and challenging sparse-angle tomography demonstrate that CR-$\ell_q$-RMM-GKS achieves significantly sharper edge preservation and higher reconstruction quality than standard non-cumulative methods, with convergence speeds competitive with or faster than non-cumulative solvers. Notably, CR-$\ell_1$-RMM-GKS outperforms CR-$\ell_2$-RMM-GKS in both edge quality and convergence speed, demonstrating that cumulative weighting and $\ell_1$ penalties are highly complementary strategies.

The remainder of the paper is organized as follows. Section~\ref{sec:background} reviews necessary background including cumulative reweighting, MM methods, and the RMM-GKS solver. Section~\ref{sec:contribution} presents the CR-$\ell_q$-RMM-GKS algorithm in detail. Section~\ref{sec:numerics} reports numerical results, and Section~\ref{sec:conclusion} concludes the paper.

\section{Background and Related Work} \label{sec:background}
This section reviews the mathematical operators and algorithmic components underlying our proposed method.

\subsection{Discrete Gradient Operator}

For a 2D image of size $n_x \times n_y$, we use the discrete gradient operator $\mat{L} \in \mathbb{R}^{2n_x n_y \times n_x n_y}$ constructed via Kronecker products:
\[
\mat{L} = \begin{bmatrix}
\mat{L}_{n_x} \otimes \mat{I}_{n_y} \\
\mat{I}_{n_x} \otimes \mat{L}_{n_y}
\end{bmatrix},
\]
where $\mat{L}_{n_x}$ and $\mat{L}_{n_y}$ are first-order finite difference operators with homogeneous Neumann boundary conditions. For 1D signals, we use $\mat{L}_n$ directly. The operator $\mat{L}\vect{x}$ approximates the image gradient at each pixel.

\subsection{\texorpdfstring{$\ell_q$}{lq}-Norm Regularization and Majorization-Minimization} \label{subsec:lq_reg}

The generalized $\ell_q$-norm regularization problem ($0 < q \leq 2$) is:
\begin{equation}
\min_{\vect{x}} J(\vect{x}) = \norm{\mat{A}\vect{x} - \vect{b}}_2^2 + \lambda \norm{\mat{L}\vect{x}}_q^q.
\label{eq:lq_problem}
\end{equation}
When $q=1$ (Total Variation-like penalty), the $\ell_1$-norm promotes sparse gradients and preserves sharp edges more effectively than $\ell_2$. For $q \leq 1$, the $\ell_q$-norm is non-differentiable, and for $q < 1$ it becomes non-convex, making direct optimization challenging.

The Majorization-Minimization (MM) principle~\cite{hunter2004, daubechies2010} addresses this by iteratively transforming the non-smooth problem into a sequence of weighted $\ell_2$ problems:
\begin{equation}
\vect{x}^{(k+1)} = \argmin_{\vect{x}} \left\{ \norm{\mat{A}\vect{x} - \vect{b}}_2^2 + \lambda \norm{\mat{W}^{(k)}\mat{L}\vect{x}}_2^2 \right\}.
\label{eq:mm_weighted_l2}
\end{equation}
Here, $\mat{W}^{(k)}$ is a diagonal weighting matrix derived from the current solution $\vect{x}^{(k)}$. Unlike the cumulative weights discussed in Section~\ref{subsec:cr}, MM weights are \emph{non-cumulative}---they depend only on the current iterate and do not accumulate edge information across iterations.

\subsection{Recycled MM-GKS Solver} \label{subsec:rmmgks}

The Majorization-Minimization Generalized Krylov Subspace (MM-GKS) framework \cite{lanza2015, okunola2025, pasha2023} efficiently solves $\ell_q$-regularized problems~\eqref{eq:lq_problem} by combining:
\begin{enumerate}
\item Majorization of the non-smooth $\ell_q$ term into weighted $\ell_2$ subproblems,
\item Projection onto a generalized Krylov subspace,
\item Automatic regularization parameter selection (e.g., discrepancy principle, L-curve) within the projected problem.
\end{enumerate}
At each iteration, MM-GKS generates weights via majorization, projects the resulting weighted problem onto a low-dimensional Krylov subspace (making the problem tractable), automatically selects the regularization parameter on the small projected system, and expands the subspace using the residual. Importantly, projection and parameter selection occur within a \emph{single iteration loop}---avoiding the nested structure of earlier projection methods.

\paragraph{Subspace Recycling}
For large-scale problems, repeated Krylov expansions can become expensive. Recycled MM-GKS (RMM-GKS) \cite{pasha2023} improves efficiency through an \emph{enlarge-compress} strategy: expand the Krylov subspace from dimension $k_{\min}$ to $k_{\max}$ to improve the solution, then compress back to $k_{\min}$ by retaining only the most informative directions. This bounds memory growth and reduces matrix-vector products. Crucially, the compressed subspace is \emph{reused} to warm-start the next outer iteration, providing computational memory between successive solves. Algorithm~\ref{alg:rmmgks} summarizes the approach.

\begin{algorithm}[htbp]
\caption{RMM-GKS: Recycled Majorization-Minimization Generalized Krylov Subspace Method}
\label{alg:rmmgks}
\begin{algorithmic}[1]
\REQUIRE $\mat{A} \in \mathbb{R}^{m \times n}$, $\mat{L} \in \mathbb{R}^{p \times n}$, $\vect{b} \in \mathbb{R}^m$, $q \in (0,2]$, $k_{\min} < k_{\max}$, $\epsilon > 0$, initial subspace $\mat{V}_0 \in \mathbb{R}^{n \times k_{\min}}$
\ENSURE Solution $\vect{x}$, parameter $\lambda$, recycled subspace $\mat{V}$

\STATE Initialize $\vect{x}^{(1,0)} = \vect{0}$, $\mat{V} \gets \mat{V}_0$, $\lambda^{(1,0)} = 0$
\STATE Compute initial MM weights: $\mat{W}^{(1,0)} = \mathrm{diag}\!\left( \bigl(|\mat{L}\vect{x}^{(1,0)}|^2 + \epsilon^2\bigr)^{(q-2)/2} \right)$

\FOR{$i = 1, 2, \ldots$ until convergence}

    \STATE \textbf{ENLARGE:} Expand subspace from $k_{\min}$ to $k_{\max}$
    \FOR{$j = 0, 1, \ldots, k_{\max} - k_{\min} - 1$}

        \STATE Compute residual gradient:
        \[
        \vect{r}^{(i,j)} = \mat{A}^T\!\bigl(\mat{A}\vect{x}^{(i,j)} - \vect{b}\bigr)
        + \lambda^{(i,j)} \mat{L}^T \bigl(\mat{W}^{(i,j)}\bigr)^2 \mat{L}\vect{x}^{(i,j)}
        \]

        \STATE Expand subspace:
        $\mat{V} \gets \bigl[\mat{V},\; \vect{r}^{(i,j)}/\|\vect{r}^{(i,j)}\|_2\bigr]$

        \STATE Select $\lambda^{(i,j)}$ via discrepancy principle~\cite{morozov1985} or L-curve~\cite{hansen1992lcurve} and solve projected problem:
        \[
        \vect{z}^{(i,j+1)} = \argmin_{\vect{z}} \left\{
        \|\mat{A}\mat{V}\vect{z} - \vect{b}\|_2^2
        + \lambda^{(i,j)}
        \|\mat{W}^{(i,j)}\mat{L}\mat{V}\vect{z}\|_2^2
        \right\}
        \]

        \STATE Update solution: $\vect{x}^{(i,j+1)} = \mat{V}\vect{z}^{(i,j+1)}$

        \STATE Update MM weights:
        $\mat{W}^{(i,j+1)} = \mathrm{diag}\!\left( \bigl(|\mat{L}\vect{x}^{(i,j+1)}|^2 + \epsilon^2\bigr)^{(q-2)/2} \right)$

    \ENDFOR

    \STATE \textbf{COMPRESS:} Retain $k_{\min}$ most informative directions via SVD (see~\cite{pasha2023}):
    $\mat{V} \gets \mat{V}\mat{M}$

    \STATE Orthogonalize current solution and append to subspace:
    \[
    \vect{v}_{\mathrm{new}}
    = \frac{\vect{x}^{(i,\,k_{\max}-k_{\min})} - \mat{V}\mat{V}^T\vect{x}^{(i,\,k_{\max}-k_{\min})}}
           {\|\vect{x}^{(i,\,k_{\max}-k_{\min})} - \mat{V}\mat{V}^T\vect{x}^{(i,\,k_{\max}-k_{\min})}\|_2},
    \qquad
    \mat{V} \gets [\mat{V},\; \vect{v}_{\mathrm{new}}]
    \]

    \IF{$\|\vect{x}^{(i,\,k_{\max}-k_{\min})} - \vect{x}^{(i-1,\,k_{\max}-k_{\min})}\|_2 \;/\;
        \|\vect{x}^{(i-1,\,k_{\max}-k_{\min})}\|_2 < \mathrm{tol}$}
        \STATE \textbf{break}
    \ENDIF

    \STATE Set $\vect{x}^{(i+1,0)} \gets \vect{x}^{(i,\,k_{\max}-k_{\min})}$,\quad
           $\lambda^{(i+1,0)} \gets \lambda^{(i,\,k_{\max}-k_{\min}-1)}$,\quad
           $\mat{W}^{(i+1,0)} \gets \mat{W}^{(i,\,k_{\max}-k_{\min})}$

\ENDFOR

\RETURN $\vect{x}^{(i,\,k_{\max}-k_{\min})}$,\quad $\lambda^{(i,\,k_{\max}-k_{\min}-1)}$,\quad $\mat{V}$
\end{algorithmic}
\end{algorithm}

\subsection{Cumulative Reweighted \texorpdfstring{$\ell_2$}{l2} Regularization} \label{subsec:cr}

Cumulative reweighting \cite{gazzola2020} enhances Tikhonov regularization by adaptively modifying spatial weights across outer iterations. At outer iteration $\ell$, it solves:
\begin{equation}
\vect{x}^{(\ell)} = \argmin_{\vect{x}} \left\{ \norm{\mat{A}\vect{x}-\vect{b}}_2^2 + \lambda_\ell \norm{\mat{D}^{(\ell)}\mat{L}\vect{x}}_2^2 \right\}.
\label{eq:crl2}
\end{equation}
Here, $\mat{D}^{(\ell)} = \mathrm{diag}(\vect{d}^{(\ell)})$ is a diagonal weighting matrix updated cumulatively via
\begin{equation}
\vect{g}^{(\ell-1)} = \frac{|\mat{D}^{(\ell-1)}\mat{L}\vect{x}^{(\ell-1)}|}{\|\mat{D}^{(\ell-1)}\mat{L}\vect{x}^{(\ell-1)}\|_\infty}, \quad
\vect{d}^{(\ell)} = \vect{d}^{(\ell-1)} \circ (1 - \vect{g}^{(\ell-1)})^{\circ s},
\label{eq:cumulative_update}
\end{equation}
where $s>0$ is a smoothness parameter, $|\cdot|$ denotes element-wise absolute value, and ${}^{\circ s}$ denotes element-wise exponentiation. The normalized weighted gradient $\vect{g}^{(\ell-1)}$ reveals newly detected edges (components near 1) and previously detected edges or smooth regions (components near 0). The multiplicative update in~\eqref{eq:cumulative_update} ensures that once an edge is detected, its associated penalty weight is progressively reduced and ``remembered'' across iterations, preventing re-smoothing and leading to robust edge preservation.

\subsection{Joint Bidiagonalization for General-Form Tikhonov 
Problems}\label{subsec:kilmer_projection}

Each subproblem in \eqref{eq:crl2} is a large-scale general-form 
Tikhonov problem: $\min_{\vect{x}} \norm{\mat{A}\vect{x}-\vect{b}}_2^2 
+ \lambda \norm{\tilde{\mat{L}}\vect{x}}_2^2$, where 
$\tilde{\mat{L}} = \mat{D}^{(\ell)}\mat{L}$. Kilmer et 
al.~\cite{kilmer2007} proposed solving such problems via a projection 
method based on joint bidiagonalization of $(\mat{A}, \tilde{\mat{L}})$. 
This iteratively constructs a small Krylov subspace $\mat{Z}_k$ onto 
which the problem is projected, reducing it to a smaller system:
\[
\min_{\vect{w}\in\mathbb{R}^k} \left\| \begin{bmatrix} 
\mat{A}\mat{Z}_k \\ \sqrt{\lambda} \tilde{\mat{L}}\mat{Z}_k \end{bmatrix} 
\vect{w} - \begin{bmatrix} \vect{b} \\ \vect{0} \end{bmatrix} 
\right\|_2^2.
\]

\paragraph{Three-Level Nested Iteration}
A key drawback of this method, especially when embedded within a 
cumulative reweighting scheme, is its deeply nested iterative structure. 
Generating each new basis vector for $\mat{Z}_k$ requires an inner 
LSQR iteration to compute orthogonal projections. Moreover, the Krylov 
subspace $\mat{Z}_k$ is typically rebuilt from scratch at each outer 
cumulative reweighting step (when $\mat{D}^{(\ell)}$ is updated). This 
leads to \emph{three nested levels of iteration}:
\begin{enumerate}
    \item outer cumulative weight update ($\mat{D}^{(\ell)}$ update),
    \item middle-level joint bidiagonalization loop (Krylov subspace 
    expansion),
    \item inner LSQR iterations for orthogonal projections.
\end{enumerate}
This deep nesting is computationally expensive and lacks ``memory'' 
between outer reweighting iterations, as subspaces are not recycled.

\section{Cumulative Reweighted \texorpdfstring{$\ell_q$}{lq} Regularization with Subspace Recycling}
\label{sec:contribution}

Two complementary strategies underpin our method: \textbf{Cumulative reweighting} (Section~\ref{subsec:cr}) excels at edge preservation by accumulating edge information across iterations, but results in expensive nested iterations. \textbf{RMM-GKS} (Section~\ref{subsec:rmmgks}) efficiently solves $\ell_q$-regularized problems with automatic parameter selection and subspace recycling, but uses non-cumulative weights that discard edge detection history.

Our contribution integrates cumulative edge tracking with the 
computational efficiency of RMM-GKS. The key insight is that RMM-GKS can naturally accommodate the changing regularization operator $\mat{D}^{(\ell)}\mat{L}$ at each outer iteration by warm-starting from the recycled subspace, eliminating nested iterations while preserving cumulative edge detection.

\subsection{Proposed Method}

We propose \emph{CR-$\ell_q$-RMM-GKS} (Cumulative Reweighted $\ell_q$ Regularization with Recycled MM-GKS), which at outer iteration $\ell = 1,2,\ldots,N_{\text{out}}$ solves

\begin{equation}
\vect{x}^{(\ell)}
=
\argmin_{\vect{x}}
\left\{
\norm{\mat{A}\vect{x} - \vect{b}}_2^2
+
\lambda_\ell
\norm{\mat{D}^{(\ell)} \mat{L}\vect{x}}_q^q
\right\},
\qquad 0<q\le2,
\label{eq:our_method_general}
\end{equation}

where the cumulative diagonal weight matrix $\mat{D}^{(\ell)} = \mathrm{diag}(\vect{d}^{(\ell)})$ is updated via~\eqref{eq:cumulative_update}. Each subproblem is solved using RMM-GKS (Algorithm~\ref{alg:rmmgks}), which: (1) applies MM to convert the $\ell_q$ problem into weighted $\ell_2$ subproblems, (2) projects onto a Krylov subspace initialized from the previous outer iteration, (3) automatically selects $\lambda_\ell$ via discrepancy principle or L-curve, and (4) returns the solution $\vect{x}^{(\ell)}$ and a compressed subspace $\mat{V}^{(\ell)}$ for recycling.

This approach offers key advantages:
\begin{enumerate}
    \item \textbf{Reduced iteration depth:} two levels (outer cumulative 
    update + inner RMM-GKS) instead of three.
    \item \textbf{Subspace recycling} across outer iterations, 
    dramatically reducing cost.
    \item \textbf{Flexibility:} supports general $\ell_q$ penalties 
    ($0 < q \le 2$), not just $\ell_2$.
    \item \textbf{Automatic parameter selection} via discrepancy 
    principle or L-curve.
    \item \textbf{Cumulative edge preservation} quality matching Gazzola et al.'s approach in \cite{gazzola2020}.
\end{enumerate}

\textbf{Remark (Convergence Properties).}
Within each outer iteration $\ell$, the cumulative weights $\mat{D}^{(\ell)}$ 
are fixed, so RMM-GKS solves a well-posed weighted $\ell_q$ problem of the 
form~\eqref{eq:our_method_general}. The MM principle guarantees that each 
inner iterate is a descent step on the majorizing surrogate, and the 
discrepancy principle provides a principled stopping criterion for the inner 
solve. Across outer iterations, the multiplicative structure 
of~\eqref{eq:cumulative_update} ensures that detected edge weights are 
monotonically non-increasing: once $\vect{d}^{(\ell)}$ is reduced, no subsequent 
update can increase it. This provides a form of one-sided monotonicity for 
the weight sequence. A full convergence theory for the outer cumulative loop 
remains open and is left for future work.
\subsection{Algorithm}

Algorithm~\ref{alg:crlq_rmmgks} presents the complete method. The outer loop updates cumulative weights to progressively identify and preserve edges, while RMM-GKS (Line~4) warm-starts from the recycled subspace $\mat{V}^{(\ell-1)}$.

\begin{algorithm}[htbp]
\caption{CR-$\ell_q$-RMM-GKS: Cumulative Reweighted $\ell_q$ Regularization}
\label{alg:crlq_rmmgks}
\begin{algorithmic}[1]
\REQUIRE $\mat{A} \in \mathbb{R}^{m \times n}$, $\mat{L} \in \mathbb{R}^{p \times n}$, $\vect{b} \in \mathbb{R}^m$, $q \in (0,2]$, $N_{\text{out}}$, $k_{\min}, k_{\max}$, $s > 0$, $\epsilon$
\ENSURE Solution $\vect{x}^{(N_{\text{out}})}$
\STATE Initialize: $\vect{x}^{(0)} = \vect{0}$, $\vect{d}^{(1)} = \mathbf{1} \in \mathbb{R}^p$, $\mat{V}^{(0)} \in \mathbb{R}^{n \times k_{\min}}$
\FOR{$\ell = 1, 2, \ldots, N_{\text{out}}$}
    \STATE $\tilde{\mat{L}} = \mathrm{diag}(\vect{d}^{(\ell)}) \mat{L}$ 
    \STATE $[\vect{x}^{(\ell)}, \lambda_\ell, \mat{V}^{(\ell)}] = \text{RMM-GKS}(\mat{A}, \tilde{\mat{L}}, \vect{b}, q, k_{\min}, k_{\max}, \epsilon, \mat{V}^{(\ell-1)})$ 
    \STATE $\vect{g}^{(\ell)} = {|\tilde{\mat{L}}\vect{x}^{(\ell)}|}/{\|\tilde{\mat{L}}\vect{x}^{(\ell)}\|_\infty}$; \quad $\vect{d}^{(\ell+1)} = \vect{d}^{(\ell)} \circ (1 - \vect{g}^{(\ell)})^{\circ s}$
\ENDFOR
\RETURN $\vect{x}^{(N_{\text{out}})}$
\end{algorithmic}
\end{algorithm}

Table~\ref{tab:iteration_comparison} compares iteration structure of edge-preserving methods. Our method achieves cumulative reweighting's edge quality while maintaining computational efficiency competitive with or better than non-cumulative approaches. The reduction from three to two iteration levels, combined with subspace recycling, provides substantial computational savings, as demonstrated in Section~\ref{sec:numerics}.

\begin{table}[htbp]
\centering
\caption{Iteration depth comparison of edge-preserving regularization 
methods. ``Levels'' counts distinct nested loops; RMM-GKS has one 
outer loop (enlarge--compress) with no separate inner solver, 
whereas Kilmer et al.'s projection requires an inner LSQR at each 
Krylov expansion step.}
\label{tab:iteration_comparison}
\begin{tabular}{lccc}
\toprule
\textbf{Method} & \textbf{Cumulative} & \textbf{Iteration} & \textbf{Subspace} \\
 & \textbf{Weighting} & \textbf{Levels} & \textbf{Recycling} \\
\midrule
Standard $\ell_q$ (RMM-GKS) & No & 1 & Yes \\
Cumulative $\ell_2$ (Gazzola et al.) & Yes & 3 & No \\
\textbf{CR-$\ell_q$-RMM-GKS (ours)} & \textbf{Yes} & \textbf{2} & \textbf{Yes} \\
\bottomrule
\end{tabular}
\end{table}

\subsection{Why CR-\texorpdfstring{$\ell_1$}{l1}-RMM-GKS Outperforms 
CR-\texorpdfstring{$\ell_2$}{l2}-RMM-GKS: Complementary Edge-Protection 
Mechanisms}\label{subsec:why_l1}

The consistently superior performance of CR-$\ell_1$-RMM-GKS over 
CR-$\ell_2$-RMM-GKS admits a clear mechanistic explanation rooted in 
how each penalty interacts with the cumulative weighting strategy.

In CR-$\ell_2$-RMM-GKS, the $\ell_2$ penalty applies uniform smoothing 
to all gradient components, and edge protection relies \emph{entirely} 
on the cumulative weight mechanism: a gradient location at index $i$ is 
protected only after the cumulative weights have reduced ${d}_i^{(\ell)}$ 
sufficiently via the multiplicative update~\eqref{eq:cumulative_update}. 
Until that detection occurs, $\ell_2$ regularization actively suppresses 
the gradient at that location --- potentially erasing soft or ambiguous 
edges before they can be identified and remembered by the cumulative 
mechanism.

In contrast, the $\ell_1$ penalty is intrinsically sparsity-promoting. 
Via the MM reweighting~\eqref{eq:mm_weighted_l2}, the effective weight 
assigned to gradient component $i$ at inner iteration $k$ is
\begin{equation}
w_i^{(k)} = \left(|\,[\mat{L}\vect{x}^{(k)}]_i|^2 + 
\epsilon^2\right)^{-1/2},
\label{eq:mm_weights_l1}
\end{equation}
which for $q=1$ assigns \emph{larger} penalties to small gradients and 
\emph{smaller} penalties to large ones. This provides a structural 
preference for sparse gradient fields that operates regardless of the 
current cumulative weight state $\vect{d}^{(\ell)}$: even before the 
cumulative weights have identified an edge, the $\ell_1$ MM weights 
already act to preserve large gradient values.

CR-$\ell_1$-RMM-GKS therefore benefits from \emph{two complementary 
edge-protection mechanisms} acting simultaneously:
\begin{enumerate}
    \item \textbf{Intrinsic sparsity promotion} from the $\ell_1$ 
    penalty, which protects large gradients at every inner iteration 
    via the MM weights~\eqref{eq:mm_weights_l1}, and
    \item \textbf{Accumulated edge memory} from the cumulative weights, 
    which progressively locks in detected edges across outer iterations 
    via~\eqref{eq:cumulative_update}.
\end{enumerate}
CR-$\ell_2$-RMM-GKS possesses only the second mechanism. This synergy 
explains two characteristic observations from Section~\ref{sec:numerics}: 
the faster convergence of CR-$\ell_1$-RMM-GKS (edges are protected from 
the first inner iteration, reducing wasted smoothing before cumulative 
detection activates) and the pronounced staircase pattern of 
CR-$\ell_2$-RMM-GKS (many outer cumulative updates are required before 
edge protection engages, producing the plateau-then-drop behaviour 
visible in Figures~\ref{fig:main_results},~\ref{fig:2d_convergence}, 
and~\ref{fig:exp3_convergence}).

\section{Numerical Results} \label{sec:numerics}

We evaluate CR-$\ell_q$-RMM-GKS on three representative inverse problems 
spanning 1D signal recovery and 2D tomographic reconstruction. The experiments 
demonstrate superior reconstruction quality over standard $\ell_2$, $\ell_1$, 
and cumulative $\ell_2$ baselines, efficient convergence, and scalability across 
problem sizes and ill-posedness regimes. All experiments highlight 
\textbf{CR-$\ell_1$-RMM-GKS}, which consistently achieves the strongest 
empirical performance.

\paragraph{Comparison Methods}
We compare against:
\textbf{$\ell_2$-RMM-GKS} (standard Tikhonov regularization~\eqref{eq:tikhonov}),
\textbf{$\ell_1$-RMM-GKS} (non-cumulative MM-GKS with $\ell_1$ penalty),
and \textbf{CR-$\ell_2$-RMM-GKS} (cumulative reweighting with $q=2$, 
analogous to Gazzola et al.~\cite{gazzola2020} but using our RMM-GKS solver).
All cumulative methods share the same RMM-GKS framework, isolating the 
effect of the $\ell_q$ penalty choice.

All experiments are implemented in Python with automatic regularization 
parameter selection via the discrepancy principle (assuming known noise level 
$\|\boldsymbol{\eta}\|_2$). Performance is measured by the relative 
reconstruction error
$
\mathrm{RRE} = \frac{\norm{\vect{x}_{\text{recon}} - 
\vect{x}_{\text{true}}}_2}{\norm{\vect{x}_{\text{true}}}_2}.
$
Inner iterations terminate when the relative change 
${\norm{\vect{x}^{(\ell,k+1)} - \vect{x}^{(\ell,k)}}_2}\,/\,
{\norm{\vect{x}^{(\ell,k)}}_2}$ falls below $10^{-5}$, 
with a total iteration budget of 600. For all experiments, we set the cumulative weighting parameter $s$ 
in~\eqref{eq:cumulative_update} to $1$, and the Krylov subspace 
dimensions $k_{\min}$ and $k_{\max}$ in Algorithm~\ref{alg:rmmgks} 
to $5$ and $25$, respectively.
\subsection{Experiment 1: 1D Signal Deblurring}

\begin{figure}[t]
    \centering
    \includegraphics[width=0.85\textwidth]{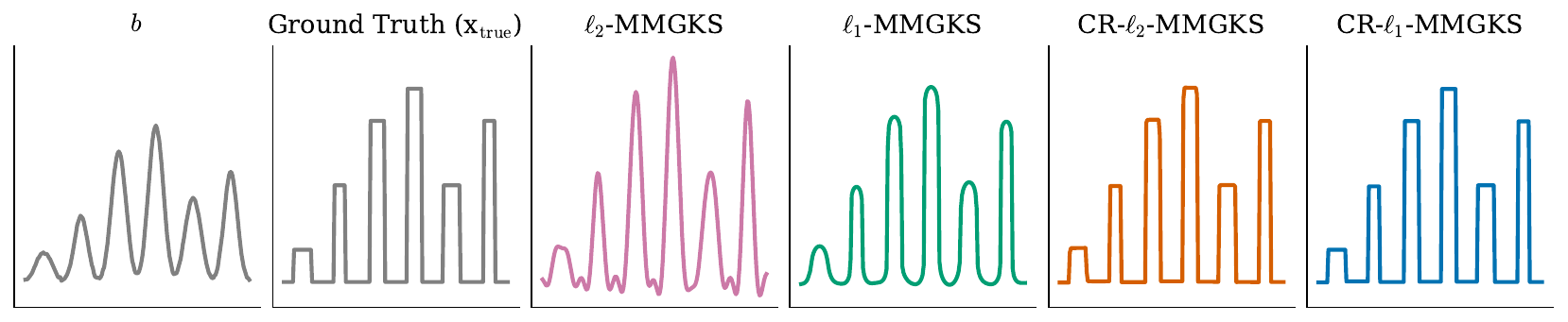}
    \caption{\textbf{Experiment 1: Final reconstructions after 600 iterations.} 
    Left to right: noisy blurred data $\vect{b}$, ground truth 
    $\vect{x}_{\text{true}}$, and reconstructions from $\ell_2$-RMM-GKS, 
    $\ell_1$-RMM-GKS, CR-$\ell_2$-RMM-GKS, and CR-$\ell_1$-RMM-GKS. 
    CR-$\ell_1$-RMM-GKS best preserves sharp jumps with minimal smoothing 
    artifacts.}
    \label{fig:snapshot}
\end{figure}

\textbf{Problem Setup.}
The goal is to reconstruct a piecewise-constant signal 
$\vect{x}_{\text{true}} \in \mathbb{R}^{200}$ containing sharp jump 
discontinuities from blurred and noisy observations. The true signal is shown 
in Figure~\ref{fig:snapshot}.

The forward operator $\mat{A} \in \mathbb{R}^{200 \times 200}$ is a Gaussian 
convolution matrix modeling a blurring point spread function, generated using 
the \textsc{TRIPs-Py} package~\cite{pasha2024trips}. Observed data 
$\vect{b} \in \mathbb{R}^{200}$ are obtained by applying $\mat{A}$ to 
$\vect{x}_{\text{true}}$ and adding white Gaussian noise $\boldsymbol{\eta}$ 
with relative noise level 
$\|\boldsymbol{\eta}\|_2 / \|\mat{A}\vect{x}_{\text{true}}\|_2 = 0.01$. 
The noisy blurred signal is shown in the leftmost column of 
Figure~\ref{fig:snapshot}. This experiment isolates each method's ability to 
recover sharp transitions while suppressing noise in smooth regions — a 
setting where $\ell_1$-type penalties have a well-known structural advantage 
over $\ell_2$.

\begin{figure}[t]
    \centering
    \begin{subfigure}{0.44\linewidth}
        \centering
        \includegraphics[width=0.95\linewidth]{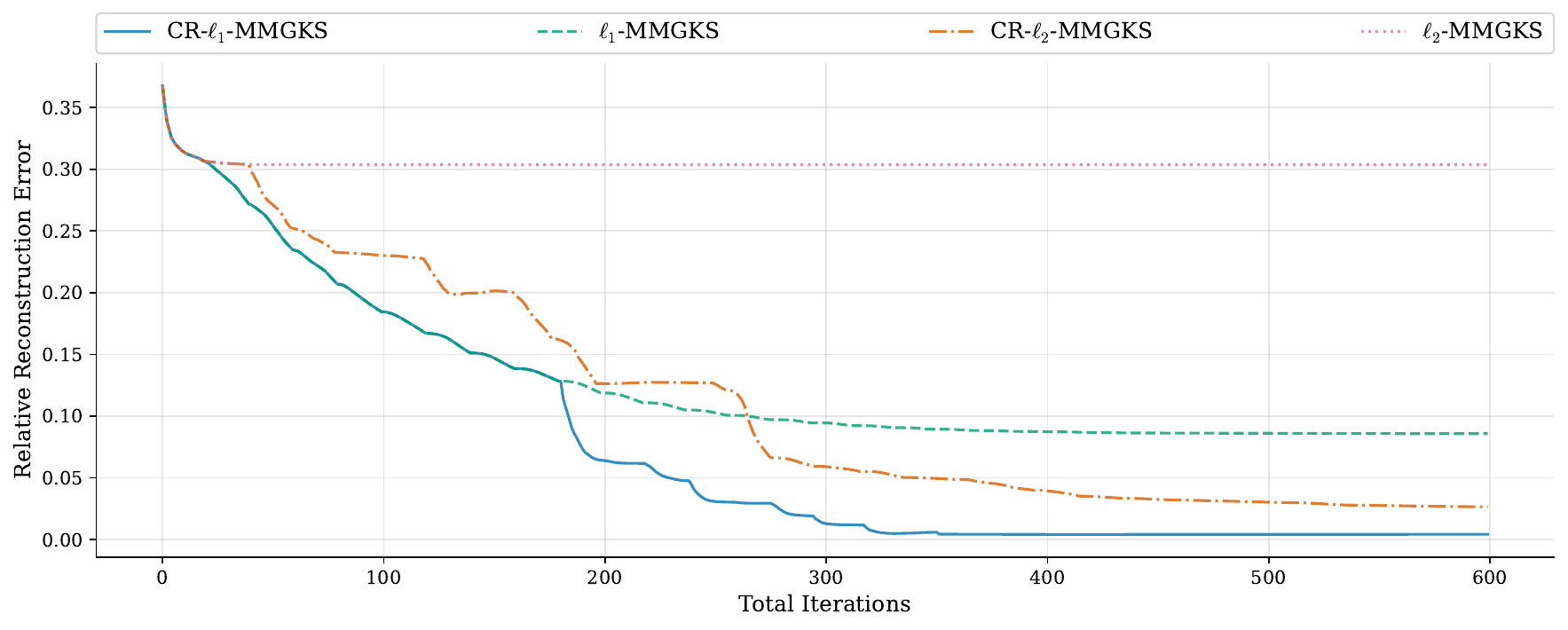}
        \caption{Relative reconstruction error versus total iterations. 
        CR-$\ell_1$-RMM-GKS (blue solid) achieves lowest error.}
    \end{subfigure}
    \hfill
    \begin{subfigure}{0.44\linewidth}
        \centering
        \includegraphics[width=0.95\linewidth]{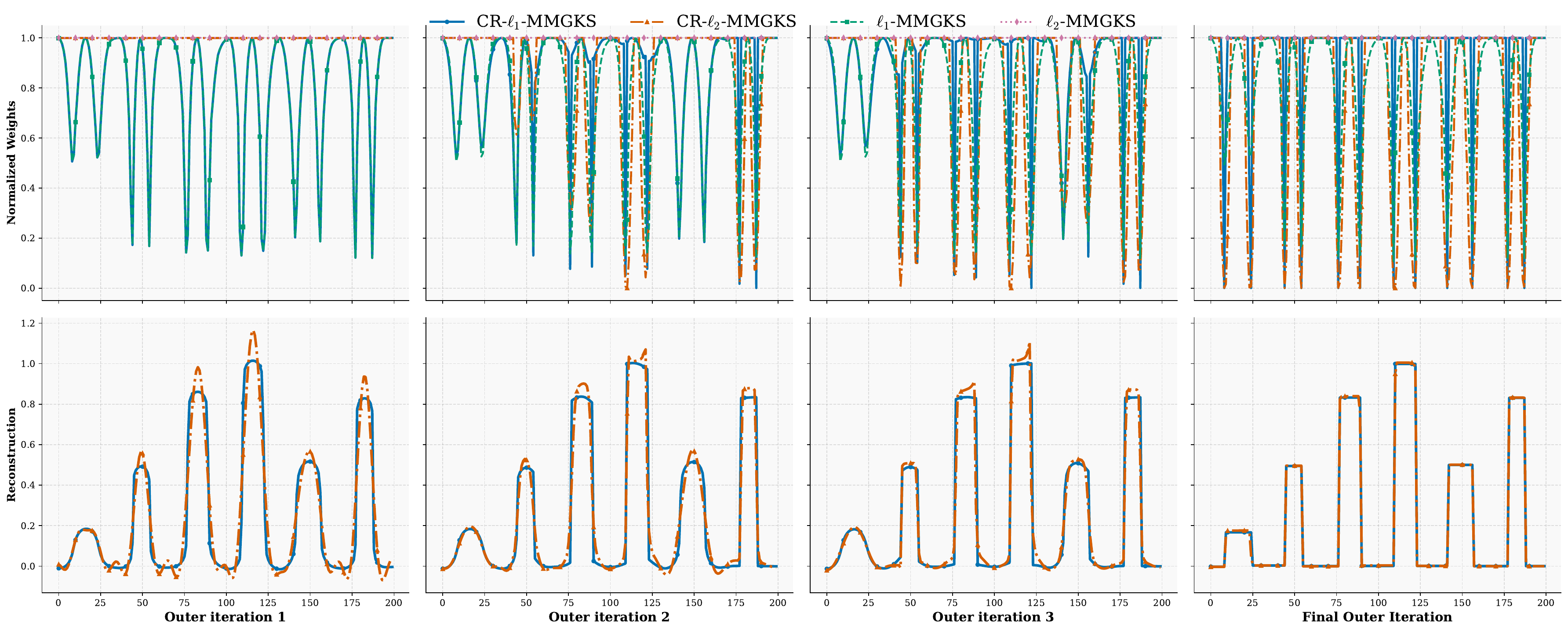}
        \caption{Evolution of cumulative weights (top row) and reconstructions 
        (bottom row) across outer iterations.}
    \end{subfigure}
    \caption{\textbf{Experiment 1: 1D deblurring.} CR-$\ell_1$-RMM-GKS 
    achieves the lowest RRE and stabilizes faster than competing methods 
    through cumulative edge memory.}
    \label{fig:main_results}
\end{figure}

\paragraph{Convergence Behavior}
Figure~\ref{fig:main_results}(a) shows RRE versus total iterations. 
CR-$\ell_1$-RMM-GKS (solid blue) decreases rapidly and achieves the lowest 
final RRE. CR-$\ell_2$-RMM-GKS (orange dash–dot) exhibits a characteristic 
\emph{staircase pattern}: each plateau corresponds to an inner RMM-GKS solve, 
and each drop to an outer cumulative weight update. This reveals that 
$\ell_2$-type reweighting requires substantially more outer iterations to reach 
comparable accuracy. Standard $\ell_1$ (green dashed) converges quickly 
initially but stabilizes at a higher error, unable to improve further without 
the memory of past iterates. Standard $\ell_2$ (red dotted) plateaus earliest, 
as expected from its tendency to oversmooth discontinuities.

\paragraph{Cumulative Weight Evolution}
Figure~\ref{fig:main_results}(b) illustrates the evolution of cumulative 
weights (top row) and corresponding reconstructions (bottom row) across outer 
iterations, providing direct insight into the edge-detection mechanism. 
Weights near zero indicate detected edges (reduced penalty); weights near one 
indicate smooth regions.

At the first outer iteration both CR-$\ell_1$- and CR-$\ell_2$-RMM-GKS begin 
with uniform weights $\vect{d}^{(1)} = \mathbf{1}$. By the second iteration, 
CR-$\ell_1$-RMM-GKS has already identified the major edges as sharp dips in 
the weight profile, whereas CR-$\ell_2$-RMM-GKS assigns more moderate 
penalties, producing visible residual smoothing at discontinuities. By the 
final iteration, CR-$\ell_1$-RMM-GKS has cumulatively locked in all edges: 
once a discontinuity is detected it remains protected from re-smoothing in 
subsequent iterations. This history-dependent mechanism is precisely what 
distinguishes cumulative from non-cumulative reweighting — $\ell_1$-RMM-GKS 
recomputes weights from the current iterate alone and thus cannot retain 
edge memory across outer iterations.

\paragraph{Snapshot Comparison}
Figure~\ref{fig:snapshot} shows final reconstructions after 600 iterations. 
CR-$\ell_1$-RMM-GKS most closely matches the ground truth, preserving all 
sharp jumps with minimal overshoot or smoothing. $\ell_2$-RMM-GKS exhibits 
the most pronounced oversmoothing. $\ell_1$-RMM-GKS and CR-$\ell_2$-RMM-GKS 
fall between these extremes, with CR-$\ell_1$-RMM-GKS achieving the smallest 
absolute reconstruction error across the entire signal.

\subsection{Experiment 2: Sparse-Angle Computed Tomography}

\begin{figure}[t]
\centering
\includegraphics[width=0.85\linewidth]{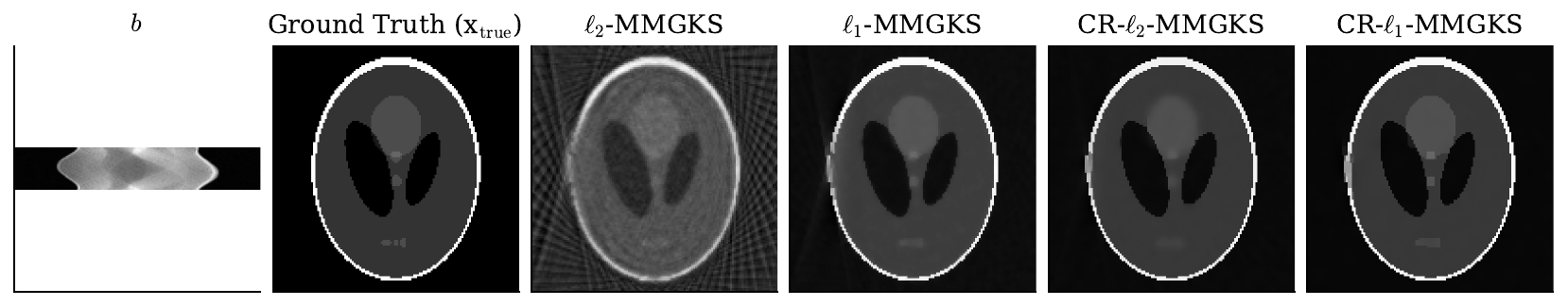}
\caption{\textbf{Experiment 2: Final reconstructions after 600 iterations.} 
Left to right: sinogram $\vect{b}$, ground truth $\vect{x}_{\text{true}}$, 
and reconstructions from $\ell_2$-RMM-GKS, $\ell_1$-RMM-GKS, 
CR-$\ell_2$-RMM-GKS, and CR-$\ell_1$-RMM-GKS. CR-$\ell_1$-RMM-GKS 
produces the sharpest tissue boundaries; $\ell_2$-RMM-GKS exhibits 
significant blurring of anatomical interfaces.}
\label{fig:2d_snapshot}
\end{figure}

\textbf{Problem Setup.}
We now move to a large-scale 2D setting. The goal is to reconstruct a simulated 
medical phantom image $\vect{x}_{\text{true}} \in \mathbb{R}^{16384}$ 
(reshaped as $128 \times 128$) from severely undersampled tomographic 
projections. The ground truth is the Shepp–Logan phantom, a standard 
benchmark in medical imaging whose elliptical regions with sharp intensity 
boundaries represent idealized tissue interfaces 
(see Figure~\ref{fig:2d_snapshot}).

The forward operator $\mat{A} \in \mathbb{R}^{5430 \times 16384}$ is a 
Radon transform matrix for fan-beam X-ray projections at only 30 
equidistant angles in $[0^\circ, 180^\circ]$~\cite{natterer2001}. Standard clinical tomography 
requires 180 or more projection angles for artifact-free reconstruction; 
our 30-angle setup is thus severely ill-posed. Observed data 
$\vect{b} \in \mathbb{R}^{5430}$ (the \emph{sinogram}) are obtained by 
applying $\mat{A}$ to the vectorized phantom and adding white Gaussian noise 
at a relative level of $1\%$. This experiment tests whether cumulative 
reweighting retains its advantage as problem scale increases and preserving 
anatomical boundaries becomes the central challenge.

\begin{figure}[t]
\centering
\includegraphics[width=0.85\linewidth]{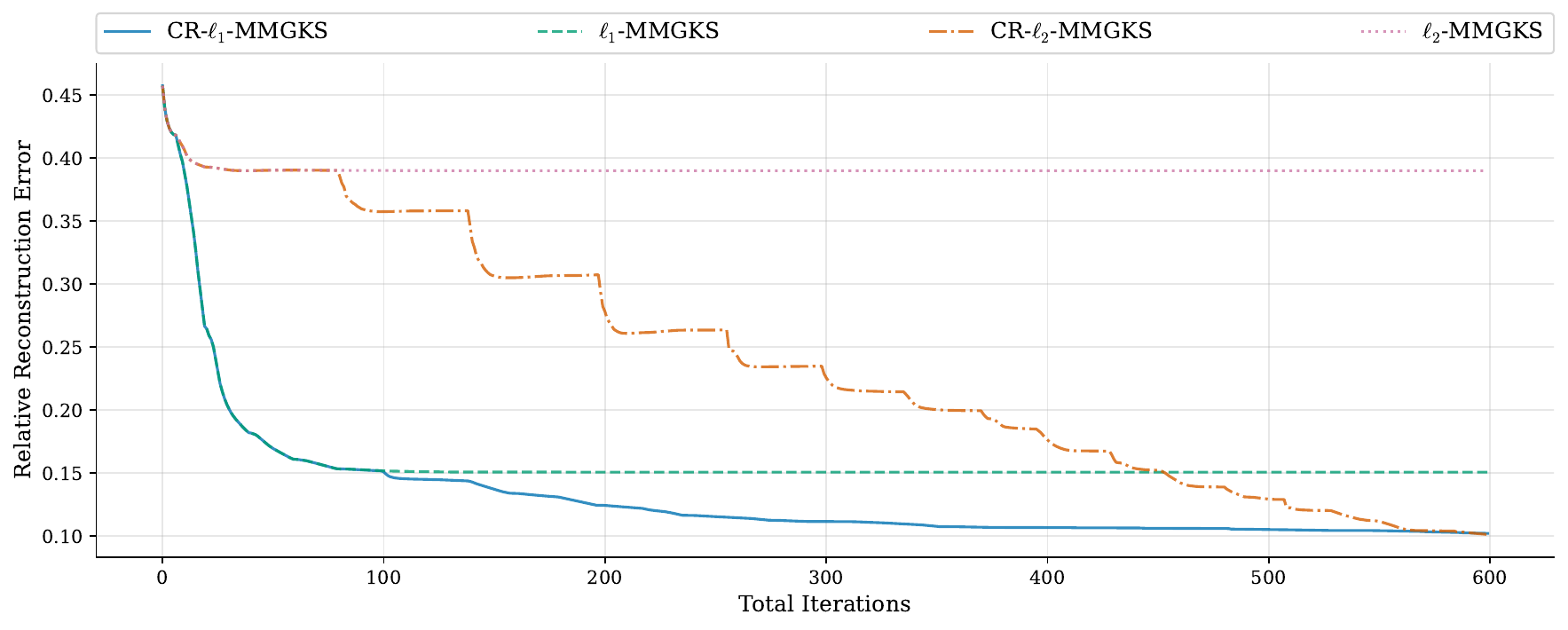}
\caption{\textbf{Experiment 2: Convergence.} CR-$\ell_1$-RMM-GKS converges 
smoothly to the lowest error. CR-$\ell_2$-RMM-GKS shows a pronounced 
staircase, requiring more outer iterations; standard $\ell_2$ flattens 
earliest due to oversmoothing.}
\label{fig:2d_convergence}
\end{figure}

\paragraph{Convergence Behavior}
Figure~\ref{fig:2d_convergence} shows RRE versus total iterations. The 
qualitative picture mirrors Experiment~1: CR-$\ell_1$-RMM-GKS achieves the 
lowest final RRE and descends most rapidly, while CR-$\ell_2$-RMM-GKS 
displays the same staircase pattern, here requiring noticeably more outer 
iterations to make progress. Standard $\ell_1$ again plateaus at an 
intermediate error, and standard $\ell_2$ flattens earliest. Importantly, 
the performance gap between cumulative and non-cumulative methods is more 
pronounced here than in the 1D case, suggesting that cumulative edge memory 
becomes increasingly valuable as problem dimension and ill-posedness grow. Although both CR methods reach the same final RRE of $0.101$ 
at the 600-iteration budget, CR-$\ell_1$-RMM-GKS achieves 
this level markedly earlier, as is evident from 
Figure~\ref{fig:2d_convergence}; the advantage of cumulative 
$\ell_1$ reweighting thus lies in convergence speed as much 
as in final accuracy for this problem.

\paragraph{Snapshot Comparison}
Figure~\ref{fig:2d_snapshot} presents final reconstructions after 600 
iterations. CR-$\ell_1$-RMM-GKS produces the sharpest boundaries between 
tissue regions, closely matching the ground truth ellipse edges. Standard 
$\ell_2$ shows significant blurring that obscures these boundaries, while 
$\ell_1$-RMM-GKS and CR-$\ell_2$-RMM-GKS recover intermediate levels of 
detail. The qualitative ranking is consistent with the RRE curves.

\subsection{Experiment 3: Limited-Angle Computed Tomography}

\begin{figure}[t]
\centering
\includegraphics[width=0.85\linewidth]{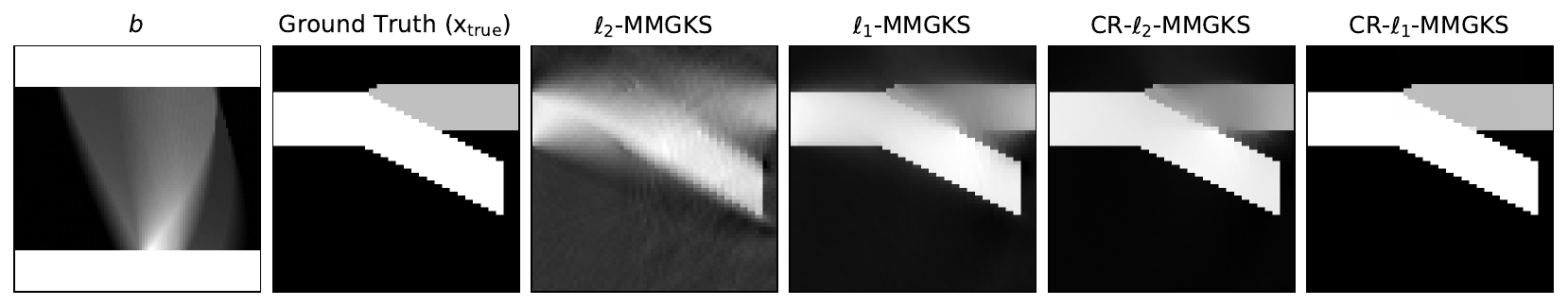}
\caption{\textbf{Experiment 3: Final reconstructions after 600 iterations.} 
Left to right: sinogram $\vect{b}$, ground truth $\vect{x}_{\text{true}}$, 
and reconstructions from $\ell_2$-RMM-GKS, $\ell_1$-RMM-GKS, 
CR-$\ell_2$-RMM-GKS, and CR-$\ell_1$-RMM-GKS. Streak artifacts from the 
missing angular range are prominent in $\ell_2$-RMM-GKS and substantially 
suppressed by CR-$\ell_1$-RMM-GKS.}
\label{fig:exp3_snapshot}
\end{figure}

\textbf{Problem Setup.}
The final experiment addresses a qualitatively different form of ill-posedness. 
Rather than reducing the number of projection angles uniformly over 
$[0^\circ, 180^\circ]$, we restrict the angular range entirely, acquiring 
projections only within a narrow cone. The goal is to reconstruct 
$\vect{x}_{\text{true}} \in \mathbb{R}^{4096}$ (reshaped as $64 \times 64$) 
from this angularly incomplete data. The ground truth is the tectonic phantom \cite{pasha2024trips}, 
a test image featuring layered structures with sharp structural boundaries
(see Figure~\ref{fig:exp3_snapshot}).

The forward operator $\mat{A} \in \mathbb{R}^{5460 \times 4096}$ encodes 
fan-beam projections at 60 equidistant angles restricted to 
$[0^\circ, 60^\circ]$ — only one third of the full angular sweep. Features 
whose dominant gradients align with unmeasured directions are invisible to 
$\mat{A}$, making recovery dependent almost entirely on the regularizer's 
structural assumptions. Observed data $\vect{b} \in \mathbb{R}^{5460}$ are 
obtained as before, with a relative noise level of $0.5\%$. This lower noise 
level reflects the focus on geometric ill-posedness rather than measurement 
corruption, allowing a cleaner comparison of how each regularizer handles 
the missing-angle problem.

\begin{figure}[t]
\centering
\includegraphics[width=0.85\linewidth]{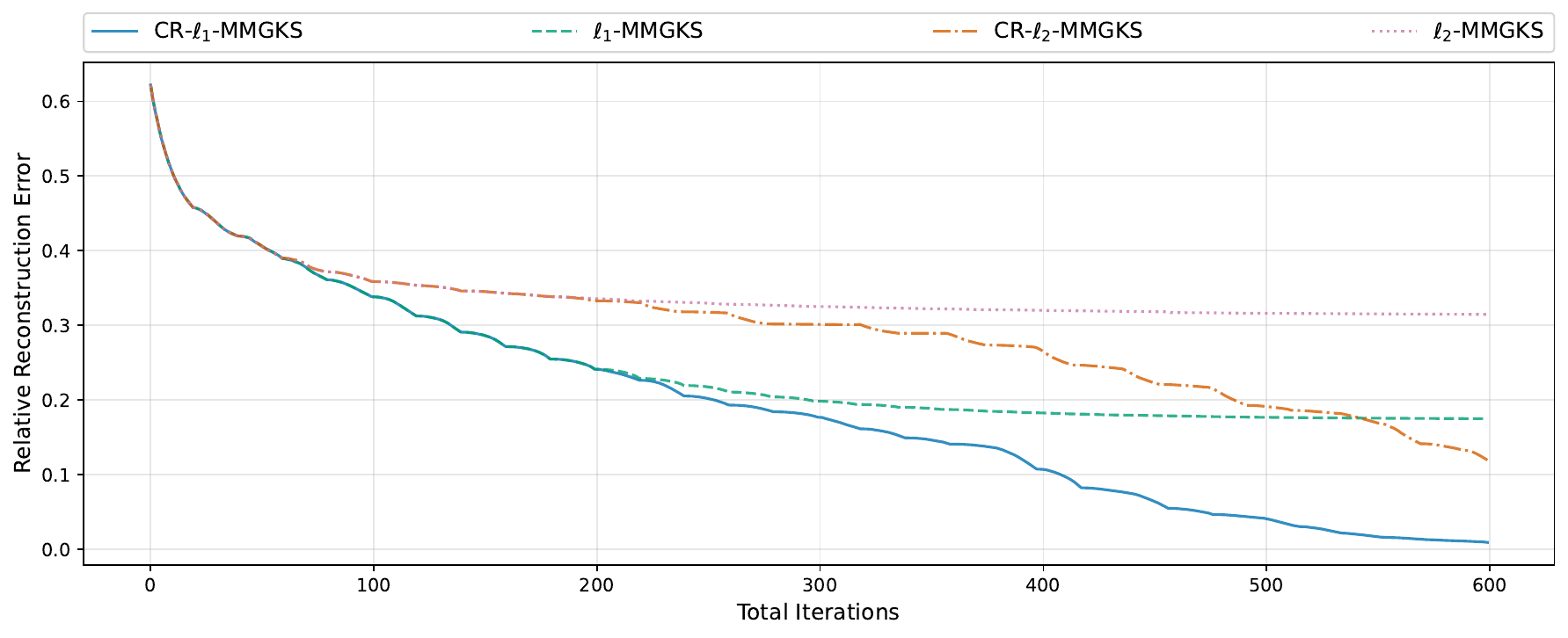}
\caption{\textbf{Experiment 3: Convergence.} CR-$\ell_1$-RMM-GKS converges 
smoothly to the lowest error. CR-$\ell_2$-RMM-GKS shows a staircase pattern; 
standard $\ell_2$ flattens at a high RRE plateau, reflecting its failure to 
resolve the limited-angle geometry.}
\label{fig:exp3_convergence}
\end{figure}

\paragraph{Convergence Behavior}
Figure~\ref{fig:exp3_convergence} shows RRE versus total iterations. 
CR-$\ell_1$-RMM-GKS converges to a markedly lower final RRE than in 
Experiment~2, consistent with the lower noise level and the strong 
sparsity-promoting structure of $\ell_1$ reweighting on piecewise-smooth 
phantoms. The staircase in CR-$\ell_2$-RMM-GKS remains visible, and the 
gap between CR-$\ell_1$- and CR-$\ell_2$-RMM-GKS is wider than in the 
sparse-angle case, indicating that cumulative $\ell_1$ reweighting is 
especially well-suited to the piecewise-flat geometry of the tectonic phantom. 
Standard $\ell_2$ flattens at a high error plateau, consistent with its 
inability to resolve the sharp directional boundaries that define this phantom.

\paragraph{Snapshot Comparison}
Figure~\ref{fig:exp3_snapshot} shows final reconstructions after 600 
iterations. The limited-angle geometry produces distinctive streak artifacts 
in the $\ell_2$-RMM-GKS reconstruction. 
CR-$\ell_1$-RMM-GKS substantially suppresses these artifacts and recovers 
the sharpest approximation to the true boundary structure. 
$\ell_1$-RMM-GKS and CR-$\ell_2$-RMM-GKS achieve partial suppression, 
falling between the two extremes and confirming that both the cumulative 
mechanism and the $\ell_1$ penalty contribute independently to artifact 
reduction.

\subsection{Summary of Findings}
Across all three experiments, the results support three consistent conclusions.
First, cumulative reweighting outperforms its non-cumulative counterpart 
regardless of the penalty choice: both CR-$\ell_1$- and CR-$\ell_2$-RMM-GKS 
achieve better edge preservation than $\ell_1$- and $\ell_2$-RMM-GKS 
respectively, confirming that retaining edge memory across outer iterations 
is a meaningful structural advantage. Second, the $\ell_1$ penalty and 
cumulative reweighting are \emph{synergistic}: CR-$\ell_1$-RMM-GKS 
consistently outperforms CR-$\ell_2$-RMM-GKS in both convergence speed and 
final reconstruction quality, demonstrating that the two mechanisms reinforce 
rather than duplicate each other. Third, this improved quality does not come 
at a computational cost: despite the added overhead of cumulative weight 
updates, CR-$\ell_1$-RMM-GKS converges in a comparable number of total 
iterations to $\ell_1$-RMM-GKS, and substantially fewer than 
CR-$\ell_2$-RMM-GKS.

\begin{table}[htbp]
\centering
\caption{Quantitative reconstruction results (RRE) after 600 
iterations. Best results per experiment in \textbf{bold}.$^{\dagger}$}
\label{tab:quantitative}
\begin{tabular}{lccc}
\toprule
\textbf{Method} & {\text{Experiment 1}} 
                & {\text{Experiment 2}}  
                & {\text{Experiment 3}} \\
\midrule
$\ell_2$-RMM-GKS        & 0.304  & 0.390 & 0.315 \\
$\ell_1$-RMM-GKS         & 0.086& 0.151 & 0.175\\
CR-$\ell_2$-RMM-GKS     &0.026 & 0.101&  0.108\\  
CR-$\ell_1$-RMM-GKS     & \textbf{0.004} & \textbf{0.101}&  \textbf{0.007} \\
\bottomrule
\end{tabular}

{\footnotesize $^{\dagger}$In Experiment~2, CR-$\ell_1$- and 
CR-$\ell_2$-RMM-GKS reach the same final RRE within the 600-iteration 
budget; CR-$\ell_1$-RMM-GKS is preferred on the basis of faster 
convergence (see Figure~\ref{fig:2d_convergence}).}
\end{table}

\section{Conclusion}\label{sec:conclusion}
We introduced CR-$\ell_q$-RMM-GKS, a framework that combines 
cumulative spatial reweighting with Krylov subspace recycling 
for edge-preserving image reconstruction. The framework 
reduces the deeply nested iteration structure of prior cumulative methods 
from three levels to two, extends cumulative weighting to general 
$\ell_q$ penalties, and retains automatic regularization parameter 
selection. Numerical experiments in 1D signal deblurring, sparse-angle tomography, 
and limited-angle tomography consistently show that CR-$\ell_1$-RMM-GKS achieves 
the lowest reconstruction error and sharpest boundaries among all 
compared methods, and that the cumulative mechanism and the $\ell_1$ 
penalty are synergistic rather than redundant.

Several limitations point toward future work. First, a formal 
convergence theory for the outer cumulative loop remains open. Second, the cumulative weighting parameter $s$ and the Krylov dimensions 
$k_{\min}, ~k_{\max}$ are currently set by hand; adaptive or 
problem-driven selection rules would improve deployability. 

\section{Acknowledgments}
M.P. gratefully acknowledges support from NSF DMS under award No. 2410699. T.O. and M.K. acknowledge support from NSF DMS No. 2410698. Any opinions, findings, conclusions, or recommendations expressed in this material are those of the authors and do not necessarily reflect the views of the National Science Foundation.

\bibliographystyle{siamplain}
\bibliography{references}

\end{document}